\documentclass[a4paper,10pt]{amsart}
\usepackage{amsfonts}
\usepackage{amsthm}
\usepackage{amssymb}
\usepackage{amsmath}
\usepackage{enumerate}
\usepackage{url}
\usepackage{color}

\begin{document}
\title{The Duffin-Schaeffer conjecture with extra divergence}

\author[C. Aistleitner]{Christoph Aistleitner}
\address{Graz University of Technology,
  Institute of Analysis and Number Theory,
  Steyrergasse 30,
  8010 Graz,
  Austria}
\email{aistleitner@math.tugraz.at}

\author[T. Lachmann]{Thomas Lachmann}
\address{Graz University of Technology,
  Institute of Analysis and Number Theory,
  Steyrergasse 30,
  8010 Graz,
  Austria}
\email{lachmann@math.tugraz.at}

\author[M. Munsch]{Marc Munsch}
\address{Graz University of Technology,
  Institute of Analysis and Number Theory,
  Steyrergasse 30,
  8010 Graz,
  Austria}
\email{munsch@math.tugraz.at}

\author[N. Technau]{Niclas Technau}
\address{Tel Aviv University,
  School of Mathematical Sciences, 
  Tel Aviv, 
  Israel}
\email{niclast@mail.tau.ac.il}

\author[A. Zafeiropoulos]{Agamemnon Zafeiropoulos}
\address{Graz University of Technology,
  Institute of Analysis and Number Theory,
  Steyrergasse 30,
  8010 Graz,
  Austria}
\email{zafeiropoulos@math.tugraz.at}

\newcommand{\mods}[1]{\,(\mathrm{mod}\,{#1})}

\begin{abstract}
The Duffin-Schaeffer conjecture is a fundamental unsolved problem in metric number theory. It asserts that for every non-negative function $\psi:~\mathbb{N} \rightarrow \mathbb{R}$ for almost all reals $x$ there are infinitely many coprime solutions $(a,n)$ to the inequality $|nx - a| < \psi(n)$, provided that the series $\sum_{n=1}^\infty \psi(n) \varphi(n) /n$ is divergent. In the present paper we prove that the conjecture is true under the ``extra divergence'' assumption that divergence of the series still holds when $\psi(n)$ is replaced by $\psi(n) / (\log n)^\varepsilon$ for some $\varepsilon > 0$. This improves a result of Beresnevich, Harman, Haynes and Velani, and solves a problem posed by Haynes, Pollington and Velani.
\end{abstract}

\maketitle

\def \l {{\lambda}}
\def \a {{\alpha}}
\def \b {{\beta}}
\def \f {{\phi}}
\def \r {{\rho}}
\def \R {{\mathbb R}}
\def \H {{\mathbb H}}
\def \N {{\mathbb N}}
\def \C {{\mathbb C}}
\def \Z {{\mathbb Z}}
\def \F {{\Phi}}
\def \Q {{\mathbb Q}}
\def \e {{\epsilon }}
\def\GL{\ensuremath{\mathop{\textrm{\normalfont GL}}}}
\def\SL{\ensuremath{\mathop{\textrm{\normalfont SL}}}}
\def\Gal{\ensuremath{\mathop{\textrm{\normalfont Gal}}}}
\def\SU{\ensuremath{\mathop{\textrm{\normalfont SU}}}}
\def\SO{\ensuremath{\mathop{\textrm{\normalfont SO}}}}

\newtheorem{prop}{Proposition}
\newtheorem{claim}{Claim}
\newtheorem{lemma}{Lemma}
\newtheorem{thm}{Theorem}
\newtheorem{defn}{Definition}
\newtheorem{conj}{Conjecture}

\theoremstyle{definition}
\newtheorem{exmp}{Example}

\theoremstyle{remark}
\newtheorem{rmk}{Remark}

\section{Introduction and statement of results}

\renewcommand{\thefootnote}{\fnsymbol{footnote}} 
\footnotetext{\emph{MSC classification:} 11J83, 11K55, 11K60.}     
\footnotetext{\emph{Keywords:} Diophantine approximation, metric number theory, Duffin-Schaeffer conjecture.}     
\renewcommand{\thefootnote}{\arabic{footnote}} 

Let $\psi: \mathbb{N} \rightarrow \mathbb{R}$ be a non-negative function. For every non-negative integer $n$ define a set $\mathcal{E}_n \subset \mathbb{R} / \mathbb{Z}$ by 
\begin{equation} \label{edef}
\mathcal{E}_n := \bigcup_{\substack{1 \leq a \leq n,\\ (a,n)=1}} \left( \frac{a - \psi(n)}{n},\frac{a+\psi(n)}{n} \right).
\end{equation}
The Lebesgue measure of $\mathcal{E}_n$ is at most $2 \psi(n) \varphi(n)/n$, where $\varphi$ denotes the Euler totient function. Thus, writing $W(\psi)$ for the set of those $x \in [0,1]$ which are contained in infinitely many sets $\mathcal{E}_n$, it follows directly from the first Borel--Cantelli lemma that $\lambda(W(\psi))=0$ whenever
\begin{equation} \label{sum}
\sum_{n=1}^\infty \frac{\psi(n) \varphi(n)}{n} < \infty.
\end{equation}
Here $\lambda$ denotes the Lebesgue measure. The corresponding divergence statement, which asserts that $\lambda(W(\psi))=1$ whenever the series in \eqref{sum} is divergent, is known as the Duffin--Schaeffer conjecture \cite{ds} and is one of the most important open problems in metric number theory. It remains unsolved since 1941.\\

The Duffin--Schaeffer conjecture is known to be true under some additional arithmetic conditions or regularity conditions on the function $\psi$. See for example \cite{har2,vaa}. In \cite{extra1} Haynes, Pollington and Velani initiated a program to establish the Duffin--Schaeffer conjecture without assuming any regularity or number-theoretic properties of $\psi$, but instead assuming a slightly stronger divergence condition. In \cite{extra1} they proved that there is a constant $c$ such that $\lambda(W(\psi))=1$, provided that
$$
\sum_{n=1}^\infty \frac{\psi(n) \varphi(n)}{n ~ e^{\left( \frac{c \log n}{\log \log n} \right)}} = \infty
$$
(throughout this paper we will understand $\log x$ as $\max(1,\log x)$, so that all appearing logarithms are positive and well-defined). In \cite{extra} Beresnevich, Harman, Haynes and Velani used a beautiful averaging argument, which is also at the core of the argument in the present paper, to show that it is sufficient to assume 
$$
\sum_{n=1}^\infty \frac{\psi(n) \varphi(n)}{n(\log n)^{\varepsilon \log \log \log n}} = \infty
$$
for some $\varepsilon > 0$. In the present paper we prove that the extra divergence factor can be reduced to $(\log n)^{\varepsilon}$ for a fixed $\varepsilon >0$. In particular this solves Problem 2 posed in \cite{extra1}, where it was asked whether the extra divergence factor $\log n$ is sufficient.

\begin{thm} \label{th1}
The Duffin--Schaeffer conjecture is true for every non-negative function $\psi: \mathbb{N} \rightarrow \mathbb{R}$ for which there is a constant $\varepsilon >0$ such that
\begin{equation} \label{div}
\sum_{n=1}^\infty \frac{\psi(n) \varphi(n)}{n(\log n)^\varepsilon} = \infty.
\end{equation}
\end{thm}

We note that by the mass transference principle of Beresnevich and Velani \cite{mtp} it is possible to deduce Hausdorff measure statements from results for Lebesgue measure, in the context of the Duffin--Schaeffer conjecture. Roughly speaking, the quantitative ``extra divergence'' result in Theorem \ref{th1} translates into a corresponding condition on the dimension function of a Hausdorff measure for the set where the Duffin--Schaeffer conjecture is true. For details we refer the reader to Section 4 of \cite{extra1}, where this connection is explained in detail. \\

\emph{Remark:} \quad While this paper was under revision, a new manuscript of Koukoulopoulos and Maynard \cite{DSCproof} appeared on the arxiv. Their paper contains a proof of the full Duffin--Schaeffer conjecture without any additional hypotheses.

\section{Proof of Theorem \ref{th1}}

In this section we give the proof of Theorem \ref{th1}. At the very end of this section, we discuss the difference between our approach and the one in \cite{extra}, and give a heuristic explanation where the improvement in the ``extra divergence'' requirement originates from.\\

Throughout the proof, we assume that $\varepsilon > 0$ is fixed. We use Vinogradov notation ``$\ll$'', where the implied constant may depend on $\varepsilon$, but not on $m,n,h$ or anything else.\\

As noted in \cite{extra}, we may assume without loss of generality that for all $n$ either $1/n \leq \psi(n) \leq 1/2$ or $\psi(n)=0$. Furthermore, by Gallagher's zero-one law \cite{gall} the measure of $W(\psi)$ can only be either 0 or 1. Thus $\lambda(W(\psi))>0$ implies $\lambda(W(\psi))=1$.\\

We will use the following version of the second Borel--Cantelli lemma (see for example \cite[Lemma 2.3]{harman}).
\begin{lemma} \label{lemmabc}
Let $\mathcal{A}_n,~n =1,2,\dots$, be events in a probability space $(\Omega,\mathcal{F},\mathbb{P})$. Let $\mathcal{A}$ be the set of $\omega \in \Omega$ which are contained in infinitely many $\mathcal{A}_n$. Assume that
$$
\sum_{n=1}^\infty \mathbb{P}(\mathcal{A}_n) = \infty.
$$
Then
$$
\mathbb{P} (\mathcal{A}) \geq \limsup_{N \to \infty} \frac{\left( \sum_{n=1}^N \mathbb{P}(\mathcal{A}_n) \right)^2}{\sum_{1 \leq m,n \leq N} \mathbb{P} (\mathcal{A}_m \cap \mathcal{A}_n)}.
$$
\end{lemma}

The following lemma of Pollington and Vaughan \cite{pv} allows us to estimate the ratio between the measure of the overlap $\mathcal{E}_m \cap \mathcal{E}_n$ and the product of the measures of $\mathcal{E}_m$ and $\mathcal{E}_n$, and is a key ingredient in \cite{extra}. In the statement of the lemma and in the sequel, we write $(m,n)$ for the greatest common divisor of two positive integers $m,n$.

\begin{lemma} \label{lemmapv}
For $m \neq n$, assume that $\lambda(\mathcal{E}_m) \lambda(\mathcal{E}_n) \neq 0$. Define
\begin{equation} \label{lemmaint}
P(m,n) = \frac{\lambda(\mathcal{E}_m \cap \mathcal{E}_n)}{\lambda(\mathcal{E}_m) \lambda(\mathcal{E}_n)}.
\end{equation}
Then
\begin{equation} \label{pdef}
P(m,n) \ll \prod_{\substack{p | \frac{mn}{(m,n)^2},\\p > D(m,n)}} \left(1 - \frac{1}{p} \right)^{-1},
\end{equation}
where the product is taken over all primes $p$ in the specified range, and where 
\begin{equation} \label{ddef}
D(m,n) = \frac{\max(n \psi(m),m \psi(n))}{(m,n)}.
\end{equation}
\end{lemma}

In view of Lemma \ref{lemmabc} it is clear that controlling $P(m,n)$ is the key to proving $\lambda(W(\psi))>0$. Following \cite{extra}, we divide the set of positive integers into blocks
\begin{equation} \label{block}
2^{4^h} \leq n < 2^{4^{h+1}}, \qquad h \geq 1,
\end{equation}
and we may assume without loss of generality that the divergence condition \eqref{div} still holds when the summation is restricted to those $n$ which are contained in a block with $h$ being even. As noted in \cite{extra}, when $m$ and $n$ are contained in different blocks, then automatically  $P(m,n) \ll 1$. Thus the real problem is that of controlling $P(m,n)$ when $m$ and $n$ are contained in the same block \eqref{block} for some $h$.\\

In the sequel, let $m,n$ be fixed, and assume that
\begin{equation} \label{*}
2^{4^h} \leq m < n < 2^{4^{h+1}}
\end{equation}
for some $h$. As in \cite{extra}, we will average the factors $P(m,n)$ over a range of downscaled versions of the sets $\mathcal{E}_m$ and $\mathcal{E}_n$. More precisely, for $k=1,2,\dots$, let $\mathcal{E}_n^{(k)}$ be defined as $\mathcal{E}_n$, but with $\psi(n) / e^k$ in place of $\psi(n)$. Correspondingly, we define
$$
P_k (m,n) = \frac{\lambda \big(\mathcal{E}_m^{(k)} \cap \mathcal{E}_n^{(k)} \big)}{ \lambda \big(\mathcal{E}_m^{(k)} \big) \lambda \big(\mathcal{E}_n^{(k)} \big)}
$$
and
$$
D_k(m,n) = \frac{\max(n \psi(m),m \psi(n))}{e^k (m,n)},
$$
and we note that for $P_k$ we have the same estimate as in \eqref{pdef}, only with $D$ replaced by $D_k$. Following \cite[Paragraph 3]{pv}, we write $m$ and $n$ in their prime factorization
$$
m = \prod_p p^{u_p}, \qquad n = \prod_p p^{v_p},
$$
and define
\begin{equation} \label{rst}
r = \prod_{\substack{p,\\u_p = v_p}} p^{u_p}, \qquad s = \prod_{\substack{p,\\u_p \neq v_p}} p^{\min(u_p,v_p)}, \qquad t = \prod_{\substack{p,\\u_p \neq v_p}} p^{\max(u_p,v_p)}.
\end{equation}
Furthermore, we set
\begin{equation} \label{dD}
\delta = \min\left( \frac{\psi(m)}{m},\frac{\psi(n)}{n} \right), \qquad \Delta = \max \left( \frac{\psi(m)}{m},\frac{\psi(n)}{n} \right).
\end{equation}
Note that
\begin{equation} \label{whatisdk}
D_k = \frac{\max(n \psi(m),m \psi(n))}{e^k (m,n)} = \frac{\textup{lcm}(m,n) \max \left( \frac{\psi(m)}{m},\frac{\psi(n)}{n} \right)}{e^k} = \Delta r t e^{-k}.
\end{equation}

Our objective is to control $P_k$ ``on average'' when summing over an appropriate range for $k$. We set 
\begin{equation} \label{K_def} 
K = K(h) = \lfloor \varepsilon h \log 4 \rfloor.
\end{equation}
Then we have the following two lemmas, which control the average size of $P_k$ for different ranges of $D_k$. We first prove the two lemmas, before we resume with the proof of Theorem \ref{th1}. The new ideas are in the proof of Lemma \ref{lemmak2}, while Lemma \ref{lemmak1} is a combination of well-known estimates.

\begin{lemma} \label{lemmak1}
Let $K$ be defined as in \eqref{K_def}. Then we have
\begin{equation} \label{comb_1}
\sum_{\substack{1 \leq k \leq K,\\  4 D_k \not\in [1,e^K)}} P_k(m,n) \ll K. 
\end{equation}
\end{lemma}

\begin{lemma} \label{lemmak2}
 Let $K$ be defined as in \eqref{K_def}. Then we have
\begin{equation*} \label{comb_2}
\sum_{\substack{1 \leq k \leq K,\\  4 D_k \in [1,e^K)}} P_k(m,n) \ll K. 
\end{equation*}
\end{lemma}

\begin{proof}[Proof of Lemma \ref{lemmak1}]
As noted in \cite{pv} and \cite{extra}, if $2 D_k \leq 1$ then $P_k(m,n) = 0$, since in this case $\mathcal{E}_m^{(k)}$ and $\mathcal{E}_n^{(k)}$ are disjoint (see the fourth displayed formula from below on p.~195 of \cite{pv}). Thus those $k$ for which $4 D_k <1$ do not contribute at all to the sum in \eqref{comb_1}. To estimate the contribution of those $k$ for which $4 D_k \geq e^K$, we note the following. By our choice of $K$ we have $e^K \gg (\log n)^\varepsilon$. Thus $4 D_k \geq e^K$ implies $D_k \gg (\log n)^\varepsilon$. As already noted in \cite{pv} and \cite{extra}, if $D_k \gg (\log n)^\varepsilon$ then we have  $P_k(m,n) \ll 1$, as a direct consequence of Lemma \ref{lemmapv} and Mertens' second theorem. This proves Lemma \ref{lemmak1}.
\end{proof}

\begin{proof}[Proof of Lemma \ref{lemmak2}]
Recall the definitions of $r,s,t$ and of $\delta$ and $\Delta$ in lines \eqref{rst} and \eqref{dD}, respectively. By the first displayed formula on page 196 of \cite{pv}, for every $k$ we have the estimate
$$
\lambda \big(\mathcal{E}_m^{(k)} \cap \mathcal{E}_n^{(k)}\big) \ll \frac{\delta}{e^k} \varphi(s) \frac{\varphi(r)^2}{r} \int_{1}^{4 \Delta r t e^{-k}} S_t(\theta)~ d\theta,
$$
where we write
$$
S_t (\theta) = \sum_{\substack{1 \leq b \leq \theta,\\ (b,t) = 1}} \frac{1}{\theta},
$$
and where we used that changing $\psi(m) \mapsto \psi(m)/e^k$ and $\psi(n) \mapsto \psi(n)/e^k$ also changes $\delta \mapsto \delta/e^k$ and $\Delta \mapsto \Delta/e^k$. Recall that by \eqref{whatisdk} we have $\Delta r t e^{-k} = D_k$. Since
$$
\lambda\big(\mathcal{E}_m^{(k)}\big) \lambda\big(\mathcal{E}_n^{(k)}\big) = \frac{\varphi(m) \varphi(n) \delta \Delta}{e^{2k}},
$$
we obtain
\begin{eqnarray*}
P_k(m,n) = \frac{\lambda \big(\mathcal{E}_m^{(k)} \cap \mathcal{E}_n^{(k)} \big)}{ \lambda \big(\mathcal{E}_m^{(k)} \big) \lambda \big(\mathcal{E}_n^{(k)} \big)} & \ll & \frac{e^k \varphi(s) \varphi(r)^2 \int_{1}^{4 D_k} S_t(\theta) ~d\theta}{\Delta r \varphi(m) \varphi(n)} \\
& = & \frac{\varphi(t) t}{\varphi(t) t} ~\frac{\varphi(s) \varphi(r)^2}{\varphi(m) \varphi(n)} ~\frac{ \int_{1}^{4 D_k} S_t(\theta) ~d\theta}{\Delta r e^{-k}} \\
& = & \frac{t}{\varphi(t)} ~\frac{\varphi(t) \varphi(s) \varphi(r)^2}{\varphi(m) \varphi(n)} ~\frac{ \int_{1}^{4 D_k} S_t(\theta) ~d\theta}{\Delta r t e^{-k}} \\
& = & \frac{t}{\varphi(t)} \frac{\int_{1}^{4 D_k} S_t (\theta) ~d\theta}{D_k},
\end{eqnarray*}
where the last line follows from $\varphi(t) \varphi(s) \varphi(r)^2 = \varphi(m) \varphi(n)$. Note that by our choice of $K$ in \eqref{K_def}, and by \eqref{*}, we have 
\begin{equation} \label{K}
e^K \ll (\log m)^{\varepsilon}, (\log n)^{\varepsilon} \ll e^K.
\end{equation}

Summing over $k$, we deduce that
\begin{equation} \label{equ_1}
\sum_{\substack{1 \leq k \leq K,\\  4 D_k \in [1,e^K)}} P_k(m,n) \ll \sum_{\substack{1 \leq k \leq K,\\  4 D_k \in [1,e^K)}}  \frac{t}{\varphi(t)} \frac{\int_{1}^{4 D_k} S_t (\theta) ~d\theta}{D_k}.
\end{equation}
We note that because of the way how $D_k$ depends on $k$, as described in \eqref{whatisdk}, there must exist a number $c \in [1,e)$ such that
$$
\left(\big\{4 D_k,~k = 1, \dots, K \big\} \cap [1,e^K) \right) \subset \{c e^j, j=0,\dots, K-1\}.
$$
Thus we have
\begin{eqnarray} \label{onther}
\sum_{\substack{1 \leq k \leq K,\\ 4 D_k \in [1,e^K)}} P_k(m,n) & \ll & \frac{t}{\varphi(t)} \sum_{j=0}^{K-1} \frac{1}{e^j}\int_{1}^{c e^j} S_t(\theta) ~d\theta.
\end{eqnarray}
For the term on the right-hand side of \eqref{onther} we have
\begin{eqnarray}
\sum_{j=0}^{K-1} \frac{1}{e^j}\int_{1}^{c e^j} S_t(\theta) ~d\theta & \ll & \sum_{j=1}^{K} \frac{1}{e^j}\int_{1}^{e^j} S_t(\theta) ~d\theta  \nonumber\\
& = & \sum_{j=1}^{K}~ \frac{1}{e^j} \sum_{\substack{1 \leq b \leq e^j,\\ (b,t) = 1}} \int_{b}^{e^j} \frac{d\theta}{\theta} \nonumber\\
& = & \sum_{j=1}^{K}~ \sum_{\substack{1 \leq b \leq e^j,\\ (b,t) = 1}} \frac{j - \log b}{e^j} \nonumber\\
& = & \sum_{\substack{1 \leq b \leq e^K,\\ (b,t) = 1}} \sum_{j = \lceil \log b \rceil}^{K} \frac{j - \log b}{e^j} \nonumber\\
& \ll & \sum_{\substack{1 \leq b \leq e^K,\\ (b,t) = 1}} \frac{1}{b} \underbrace{\sum_{i=1}^{\infty} \frac{i}{e^i}}_{\ll 1} \nonumber\\
& \ll & \sum_{\substack{1 \leq b \leq e^K,\\ (b,t) = 1}} \frac{1}{b}. \label{sieve}
\end{eqnarray}
The sum in \eqref{sieve} can be estimated using a sieve with logarithmic weights. Following the lines of \cite[Lemma 2.1]{gkm}, we have
\begin{eqnarray}
\sum_{\substack{1 \leq b \leq e^K,\\(b,t)=1}} \frac{1}{b} & = & \sum_{\substack{1 \leq b \leq e^K,\\ p | b \implies p \nmid t}} \frac{1}{b} \nonumber\\
& \leq & \prod_{\substack{p \leq e^K,\\ p \nmid t}} \left( 1 - \frac{1}{p} \right)^{-1} \nonumber\\
& = & \left(\prod_{p \leq e^K} \left( 1 - \frac{1}{p} \right)^{-1}\right)  \left(\prod_{\substack{p \leq e^K,\\ p \mid t}} \left( 1 - \frac{1}{p} \right)\right). \label{lastline}
\end{eqnarray}
For the first product in \eqref{lastline} by Mertens' theorem we have
$$
\prod_{p \leq e^K} \left( 1 - \frac{1}{p} \right)^{-1} \ll K.
$$
For the second product we have
\begin{eqnarray*}
\prod_{\substack{p \leq e^K,\\ p \mid t}} \left( 1 - \frac{1}{p} \right) & = & \frac{\varphi(t)}{t} \underbrace{\prod_{\substack{p > e^K,\\ p \mid t}} \left( 1 - \frac{1}{p} \right)^{-1}}_{\ll 1} \ll \frac{\varphi(t)}{t},
\end{eqnarray*}
where Mertens' theorem and \eqref{K} were used. Finally, inserting these bounds into \eqref{sieve} and combining this with \eqref{onther} establishes Lemma \ref{lemmak2}.
\end{proof}

\begin{proof}[Conclusion of the proof of Theorem \ref{th1}:] 
By a combination of Lemma \ref{lemmak1} and Lemma \ref{lemmak2} we immediately obtain
\begin{equation} \label{p_est}
\sum_{k=1}^K P_k(m,n) \ll K.
\end{equation}
By the definition of $P_k(m,n)$ we have
\begin{eqnarray*}
\sum_{k=1}^K P_k(m,n) & = & \sum_{k=1}^K \frac{\lambda\big(\mathcal{E}_m^{(k)} \cap \mathcal{E}_n^{(k)}\big)}{\lambda \big(\mathcal{E}_m^{(k)} \big) \lambda \big(\mathcal{E}_n^{(k)} \big)} \\
& = & \sum_{k=1}^K \frac{e^{2k} \lambda\big(\mathcal{E}_m^{(k)} \cap \mathcal{E}_n^{(k)}\big)}{\lambda \big(\mathcal{E}_m \big) \lambda \big(\mathcal{E}_n \big)},
\end{eqnarray*}
and consequently \eqref{p_est} implies that
$$
\sum_{k=1}^K e^{2k} \lambda\big(\mathcal{E}_m^{(k)} \cap \mathcal{E}_n^{(k)}\big) \ll K \lambda \big(\mathcal{E}_m \big) \lambda \big(\mathcal{E}_n \big).
$$
The implied constant is independent of $m$ and $n$, and thus summing over $m$ and $n$ yields 
$$
\sum_{k=1}^K ~\sum_{2^{4^h} \leq m < n < 2^{4^{h+1}}} e^{2k} \lambda\big(\mathcal{E}_m^{(k)} \cap \mathcal{E}_n^{(k)}\big) \ll K \sum_{2^{4^h} \leq m < n < 2^{4^{h+1}}} \lambda \big(\mathcal{E}_m \big) \lambda \big(\mathcal{E}_n \big).
$$
Accordingly, there is at least one choice of $k=k(h)$ in the range $\{1, \dots, K\}$ such that
\begin{equation*}
\sum_{2^{4^h} \leq m < n < 2^{4^{h+1}}} e^{2k} \lambda\big(\mathcal{E}_m^{(k)} \cap \mathcal{E}_n^{(k)}\big) \ll \sum_{2^{4^h} \leq m < n < 2^{4^{h+1}}} \lambda \big(\mathcal{E}_m \big) \lambda \big(\mathcal{E}_n \big),
\end{equation*}
or, equivalently, such that
\begin{equation} \label{indep}
\sum_{2^{4^h} \leq m < n < 2^{4^{h+1}}} \lambda\big(\mathcal{E}_m^{(k)} \cap \mathcal{E}_n^{(k)}\big) \ll \sum_{2^{4^h} \leq m < n < 2^{4^{h+1}}} \lambda \big(\mathcal{E}_m^{(k)} \big) \lambda \big(\mathcal{E}_n^{(k)} \big),
\end{equation}
where the implied constant does not depend on $h$. We replace the original function $\psi(n)$ by a function $\psi^*(n)$, where
$$
\psi^*(n) = \left\{ \begin{array}{ll} 0 & \text{when $n$ is not in $\big[2^{4^h}, 2^{4^{h+1}} \big)$ for some even $h$},\\ \psi(n) e^{-k(h)} & \text{when $n$ is in $\big[2^{4^h}, 2^{4^{h+1}} \big)$ for some even $h$,} \end{array} \right.
$$
and write $\mathcal{E}_n^*, ~n \geq 1,$ for the corresponding sets, which are defined as in \eqref{edef} but with $\psi^*$ in place of $\psi$. By \eqref{K} we have
$$
\psi^*(n) \gg \frac{\psi(n)}{(\log n)^\varepsilon}.
$$
Thus the extra divergence condition in the assumptions of Theorem \ref{th1} guarantees that
$$
\sum_{n=1}^\infty \lambda(\mathcal{E}_n^*) = \infty,
$$
while \eqref{indep} guarantees that
$$
\sum_{1 \leq m, n \leq N} \lambda \big(\mathcal{E}_m^* \cap \mathcal{E}_n^* \big) \ll \sum_{1 \leq m,n \leq N} \lambda(\mathcal{E}_m^*) \lambda(\mathcal{E}_n^*)
$$
(recall here that $\lambda(\mathcal{E}_m^* \cap \mathcal{E}_n^*) \ll \lambda(\mathcal{E}_m^*) \lambda(\mathcal{E}_n^*)$ holds automatically when $m$ and $n$ are not contained in the same block for some $h$). Thus by Lemma \ref{lemmabc} we have $\lambda(W(\psi^*)) > 0$, and since $\mathcal{E}_n^* \subset \mathcal{E}_n$ we also have $\lambda(W(\psi))>0$. By Gallagher's zero-one law, positive measure of $W(\psi)$ implies full measure. Thus $\lambda(W(\psi))=1$, which proves the theorem.
\end{proof}

We end with a comparison of our proof and the one in \cite{extra}. At the core of the argument in \cite{extra} was the observation that
\begin{eqnarray}
\sum_{k=1}^K P_k(m,n) & \ll & \sum_{k=1}^K~ \prod_{\substack{p | \frac{mn}{(m,n)^2},\\p > e^k}} \left(1 - \frac{1}{p} \right)^{-1} \nonumber\\
& \ll & \sum_{k=1}^K \frac{\log \log n}{k} \nonumber\\
& \ll & (\log K)(\log \log n), \label{endup}
\end{eqnarray}
where the product in the first line was directly estimated using Mertens' second theorem. Thus when $K \gg (\log \log n)(\log \log \log n)$ we have $\sum_{k=1}^K P_k(m,n) \ll K$, and accordingly there is at least one value of $k$ in this range for which $P_k(m,n) \ll 1$. Extending this argument over a range of pairs of integers $m$ and $n$ instead of assuming that $m,n$ are fixed allows to deduce the desired ``extra divergence'' result, with the requirement of potentially having to divide $\psi(n)$ by a factor as large as $e^K \leq e^{\varepsilon (\log \log n)(\log \log \log n)}$ for all $n$ and still keeping the divergence of the series of measures.\\

In our proof we followed roughly  the same plan. However, instead of taking the estimate 
$$
\prod_{\substack{p | \frac{mn}{(m,n)^2},\\p > e^k}} \left(1 - \frac{1}{p} \right)^{-1} \ll \frac{\log \log n}{k}
$$ 
for granted and then averaging over different reduction factors $e^k$, we took the averaging procedure into the sieve argument which is at the core of the proof of Lemma \ref{lemmapv}. To see why an improvement is possible by proceeding this way, we note that to obtain the overlap estimate in Lemma \ref{lemmapv} it is necessary to give upper bounds for sums 
$$
\sum_{\substack{1 \leq b \leq \theta,\\(b,t)=1}} 1,
$$
where $\theta \ll \log t$ essentially is $D(m,n)$ from \eqref{ddef}, and where $t$ is the quotient $\frac{mn}{(m,n)^2}$ which appears in \eqref{pdef}. The purpose of Lemma \ref{lemmapv} is to compare the size of this sifted sum to $\theta \varphi(t) / t$. To obtain Lemma \ref{lemmapv} one applies the classical sieve bound
\begin{equation} \label{prod_r}
\sum_{\substack{1 \leq b \leq \theta,\\(b,t)=1}} 1 \ll \theta \prod_{\substack{p \mid t,\\p \leq \theta}} \left(1 - \frac{1}{p} \right) = \theta \frac{\varphi(t)}{t} \prod_{\substack{p \mid t,\\p > \theta}} \left(1 - \frac{1}{p} \right)^{-1},
\end{equation}
and the product on the right is the one which also appears in \eqref{pdef}. This sieve bound gives optimal results for some constellations of parameters, but in our proof we used the fact that we were actually averaging over different values of $k$ (which determine $\theta$) to save some factors. We exhibit two extremal cases showing this phenomenon. The factor $P(m,n)$ can only be large when the product on the right of \eqref{prod_r} is large. However, this product can only be large if a very large proportion of small primes divides $t$. Assume on the contrary that \emph{no} small prime divides $t$. Then the sieve inequality in \eqref{prod_r} is optimal since it actually is an equality, since on both sides we have exactly $\theta$. However, in this case the product on the right is extremely small and cannot cause problems. As a second extremal case, assume that \emph{all} small primes divide $t$. Then the product on the right is very large, but the sieve bound is not sharp, since in the sum on the left the only number we count is the number 1 (no other small number is coprime to $t$). So there is a trade-off between the way a large proportion of primes dividing $t$ is able to increase the value of the product on the right of \eqref{prod_r}, but at the same time reduces the quality of the sieve bound, and this can be exploited when summing over different values of $\theta$ in \eqref{prod_r}. One might expect that this should be a very subtle relationship, and in general this is indeed the case (cf. \cite[Proposition 2.6]{gkm}, where this phenomenon is addressed). However, quite surprisingly, it turned out that in our particular setting it was possible to exploit this phenomenon by a relatively simple calculation.

\section*{Acknowledgements}

CA is supported by the Austrian Science Fund (FWF), projects F-5512, I-3466 and Y-901. NT is supported by FWF projects W-1230 and Y-901. TL and MM are supported by FWF project Y-901. AZ is supported by FWF projects F-5512 and Y-901. This paper is an outcome of the Fufu Seminar.


\end{document}